# The Legendre Spectral-Collocation method for a class of fractional integral equations

A. Yousefi, S. Javadi, E. Babolian

Department of Computer Science, Faculty of Mathematical Sciences and Computer, Kharazmi University, Tehran, Iran

asadyosefi@gmail.com, javadi@khu.ac.ir, babolian@khu.ac.ir

**Abstract**

In this paper, we consider spectral-collocation method base on Legendre-Gauss-Lobatto point. We present a computational method for solving a class of fractional integral equation of the second kind. Then based on Legendre-Gauss-Lobatto point and using, we derive a system of algebraic equations. The method is illustrated by applications and the results obtained are compared with the exact solutions in open literature. The obtained numerical results show that our proposed method is efficient and accurate for fractional integral equations of second kind. In addition, we prove that the error of the approximate solution decay exponentially in $L^2$ norm.

**Keyword.** Legendre Spectral-Collocation method, Fractional integral equation, Fractional order.

## 1. Introduction

Fractional calculus has gained importance and popularity during the past three decades or so, due to mainly its demonstrated applications in numerous seemingly diverse fields of science and engineering. The advantage of fractional calculus becomes apparent in modeling mechanical and electrical properties of real material, as well as in the description of many fractional fields (see for details [1-2]).

In recent years, there has been many studies of fractional integral and differential equations (for instance, see [3-13]). Spectral methods are an emerging area in the field of applied science and engineering. These methods provide a computational approach that has achieved substantial popularity over the last three decades. They have been applied successfully to numerical simulations of many problems in fractional calculus ([14-19]).

In this paper, we are concerned with the numerical solutions of the following equation:
$$y(x) = a(x)\, I^\alpha\{b(x)y(x)\} + f(x), \qquad (1)$$
where $I^\alpha$ is the fractional integral of order $0 < \alpha < 1$, $a, b$ and $f$ are the continuous real functions on $[0, T]$, $y$ is the unknown function that must be determined. Existence and uniqueness of the solution of the Eq.(1), in additional to some analytical properties and important inequalities, are investigated [20]. The authors in [20] applied the product trapezoidal method based on Picard iteration for approximating solution at the mesh point.

To the best of our knowledge, fractional derivative and integral are global, i.e. they are defined by an interval over the whole interval $[0, T]$, and therefore global method, such as spectral methods, are better suited for this equation. We introduce a collocation method based on the Gauss-Legendre interpolation for solving Eq. (1). Inspired by the work of [21], we extend the approach to Eq. (1) and provide a rigorous convergence analysis for the Legendre spectral-collocation method. We show that an approximate solution is convergent in $L^2$ and $L^\infty$ norms. The results obtained in this paper may be considered as a generalization of the result in [20].

The structure of this paper is as follows: In section 2, some necessary definitions and mathematical tools of the fractional calculus which are required for our subsequent developments are introduced. In section 3, the Legendre spectral-collocation method of the fractional integral equation of second kind is obtained. After this section, we discuss about convergence analysis and then, some numerical experiments are presented in Section 5 until efficiency of Legendre spectral-collocation method. Some conclusions are given in section 6.

## 2. Basic Definitions and Notations

For the concept of fractional integral, we will adopt Riemann-Liouville definition. In whole of this paper, integral operator will be Riemann-Liouville fractional integral operator.

**Definition 1.2** A real function $f$ on $[0, T]$ is said to be in the space $C_\mu$, $\mu \in \mathbb{R}$, if there exists a real number $p(> \mu)$, such that $f(x) = x^p f_1(x)$, where $f_1 \in C[o, T]$, and it is said to be in the space $C_\mu^m$ iff $f^{(m)} \in C_\mu$, $m \in \mathbb{N}$.

**Definition 2.2** The Riemann-Liouville fractional integral operator of order $\alpha \geq 0$, of a function $f \in C_\mu$, $\mu \geq -1$, is defined as

$$I^\alpha f(x) = \frac{1}{\Gamma(\alpha)} \int_0^x (x-t)^{\alpha-1} f(t) dt, \qquad \alpha > 0, \qquad x \in [0, T], \tag{2}$$

$$I^0 f(x) = f(x),$$

where $\Gamma(.)$ is the Gamma function where

$$\Gamma(\alpha) = \int_0^\infty t^{\alpha-1} e^{-t} dt.$$

Some properties of the Riemann-Liouville fractional integral are

$$(i)\ I^\alpha I^\beta = I^{\alpha+\beta}, \quad (ii)\ I^\alpha I^\beta = I^\beta I^\alpha,$$

$$(iii)\ I^\alpha (x-a)^\nu = \frac{\Gamma(\nu+1)}{\Gamma(\alpha+\nu+1)} (x-a)^{\alpha+\nu},$$

where $\alpha, \beta \geq 0$ and $\nu > -1$.

Let $\omega^{q_1,q_2}(x) = (1-x)^{q_1}(1+x)^{q_2}$ be a weight function in the usual sense, for $q_1, q_2 > -1$. The Jacobi polynomials $\{J_n^{q_1,q_2}\}_{n=0}^\infty$ are orthogonal with respect to $\omega^{q_1,q_2}$ over interval $(-1,1)$. The set of Jacobi polynomials is a complete $L^2_{\omega^{q_1,q_2}}(-1,1)$ −orthogonal system, namely [17]

$$\int_{-1}^1 J_n^{q_1,q_2} J_m^{q_1,q_2} \omega^{q_1,q_2}(x)\, dx = \gamma_n^{q_1,q_2} \delta_{nm}, \tag{3}$$

where $\delta_{nm}$ is the Kronecker function, and

$$\gamma_n^{q_1,q_2} = \begin{cases} \dfrac{2^{q_1+q_2+1}\,\Gamma(q_1+1)\Gamma(q_2+1)}{\Gamma(q_1+q_2+2)}, & n = 0, \\ \dfrac{2^{q_1+q_2+1}\,\Gamma(n+q_1+1)\Gamma(n+q_2+1)}{(2n+q_1+q_2+1)n!\,\Gamma(n+q_1+q_2+1)}, & n \geq 1. \end{cases}$$

In particular, we find

$$J_0^{q_1,q_2}(x) = 1, \quad J_1^{q_1,q_2}(x) = \frac{1}{2}(q_1+q_2+2)x + \frac{1}{2}(q_1-q_2).$$

Let $S_N[-1,1]$ the set of all polynomials of degree $\leq N$ ($N \geq 0$). Thus, if we denote by $\{x_j^{q_1,q_2}, \omega_j^{q_1,q_2}\}_{j=0}^{N}$ the nodes and the corresponding Christoffel numbers of the standard Jacobi-Gauss interpolation on $\Lambda = (-1,1)$, then we have

$$\int_{-1}^{1} v(x)\omega_j^{q_1,q_2}(x)dx = \sum_{j=0}^{N} v(x_j^{q_1,q_2})\omega_j^{q_1,q_2}, \quad \forall v \in S_{2N-1}[-1,1]. \tag{4}$$

For any $u \in C(\Lambda)$, we denote by $p_{x,N}^{q_1,q_2} : C(\Lambda) \to S_N$ the Jacobi-Gauss interpolation operator, such that

$$p_{x,N}^{q_1,q_2} u(x_j^{q_1,q_2}) = u(x_j^{q_1,q_2}), \quad 0 \leq j \leq N. \tag{5}$$

It is clear that

$$p_{x,N}^{q_1,q_2} u(x) = \sum_{j=0}^{N} u_j^{q_1,q_2} J_j^{q_1,q_2}(x), \tag{6}$$

where

$$u_j^{q_1,q_2} = \frac{1}{\gamma_j^{q_1,q_2}} \sum_{j=0}^{N} u(x_j^{q_1,q_2}) J_j^{q_1,q_2}(x_j^{q_1,q_2}) \omega_j^{q_1,q_2}.$$

In special case, if $q_1 = q_2 = 0$, then the Jacobi polynomial is reduced to the Legendre polynomial. Thus, we recall that the Legendre polynomials [17-18] $\{L_i(x); i = 0,1, \dots\}$ are defined on the $\Lambda$ with the following recurrence formula:

$$L_{i+1}(x) = \frac{2i+1}{i+1} x\, L_i(x) - \frac{i}{i+1} L_{i-1}(x), \quad i = 1,2,\dots,$$

with $L_0(x) = 1$ and $L_1(x) = x$. By using the Rodrigues formula $L_i(x) = \frac{(-1)^i}{2^i\, i!} D^k\big((1-x^2)^i\big)$, the Legendre polynomial has the expansion as follow

$$L_i(x) = \frac{1}{2^i} \sum_{k=0}^{\left[\frac{i}{2}\right]} (-1)^k \frac{(2i-2k)!}{2^k k!\,(i-k)!\,(i-2k)!} x^{i-2k}. \tag{7}$$

The set of $\{L_i(x); i = 0,1, \dots\}$ is a complete orthogonal system in $L^2(\Lambda)$, and we have

$$(L_j, L_k) = \int_{-1}^{1} L_j(x) L_k(x) \, dx = h_j \delta_{jk}, \tag{8}$$

where $\delta_{jk}$ is the Kronecker delta symbol and $h_j = \frac{2}{2j+1}$. Thus, for any $g \in L^2(\Lambda)$, we obtain

$$g(x) = \sum_{i=0}^{\infty} a_i L_i(x), \quad a_i = \frac{1}{h_i} \int_{-1}^{1} g(x) L_i(x) \, dx. \tag{9}$$

For simplicity, we let $x_j = x_j^{0,0}$, $\omega_j = \omega_j^{0,0}$ and $p_{x,N} = p_{Nx,}^{0,0}$. Thus, according to Eq.(4), we have

$$\int_{-1}^{1} v(x) dx = \sum_{j=0}^{N} v(x_j) \omega_j, \quad \forall v \in S_{2N-1}[-1,1], \tag{10}$$

where $x_j$ $(0 \leq j \leq N)$, i.e. are the zero of $(1-x^2) L'_N(x)$, and $\omega_j$ $(0 \leq j \leq N)$, i.e. $\omega_j = \frac{2}{N(N+1)\left(L_N(x_j)\right)^2}$, are denoted to the nodes and Christoffel numbers of Legendre-Gauss-Lobatto interpolation on the classical interval $(-1,1)$, repectively. The norm and discrete inner product in $L^2(\Lambda)$ are defined as

$$(u,v)_N = \sum_{j=0}^{N} u(x_j) v(x_j) \omega_j, \quad \|u\|_N = (u,u)^{\frac{1}{2}}. \tag{11}$$

In the next section, we use some special techniques for converting Eq. (1) to a problem can be solved by the standard Legendre spectral collocation scheme. These topics will be presented in the next section.

## 3. The Spectral Collocation Method

In this section, we shall propose a Legendre spectral collocation method for solving the Eq. (1). The Legendre spectral-collocation method takes advantage of both the Legendre polynomials and Legendre-Gauss-Lobatto interpolation nodes. The main idea is to use the suggested method consists of discretizing the fractional integral equation to obtain a system of algebraic equations with unknown coefficients. As a result, the solution of the obtained system can be easily solved.

Again, let us to consider Eq. (1) by applying **definitions 2.2** as follow

$$y(x) = \frac{a(x)}{\Gamma(\alpha)} \int_0^x b(s)(x-s)^{\alpha-1} y(s) \, ds + f(x), \quad x \in I = [0,T].$$

Next, in enjambment we let $\Lambda = [-1,1]$. For ease of analysis, we transfer the problem (1) to an equivalent problem defined in $\Lambda$. More specifically, we apply the change of variable $t = \frac{2x}{T} - 1$ or $x = \frac{T}{2}(t+1)$, such that $t \in \Lambda$. Then, we have

$$y\left(\frac{T}{2}(t+1)\right) = \frac{a\left(\frac{T}{2}(t+1)\right)}{\Gamma(\alpha)} \int_0^{\frac{T}{2}(t+1)} b(s) \left(\frac{T}{2}(t+1) - s\right)^{\alpha-1} y(s)\, ds$$
$$+ f\left(\frac{T}{2}(t+1)\right), \quad t \in \Lambda. \tag{12}$$

Moreover, to transfer the integral interval $\left(0, \frac{T}{2}(t+1)\right)$ to $(-1, t)$, we make the transformation $s = \frac{T}{2}(\mu + 1)$. Then Eq. (12) can be written

$$y\left(\frac{T}{2}(t+1)\right) = \left(\frac{T}{2}\right)^\alpha \frac{a\left(\frac{T}{2}(t+1)\right)}{\Gamma(\alpha)} \int_{-1}^t b\left(\frac{T}{2}(\mu+1)\right)(t-\mu)^{\alpha-1} y\left(\frac{T}{2}(\mu+1)\right) d\mu$$
$$+ f\left(\frac{T}{2}(t+1)\right), \quad t \in \Lambda. \tag{13}$$

Further, if we let

$$Y(t) = y\left(\frac{T}{2}(t+1)\right), \quad A(t) = a\left(\frac{T}{2}(t+1)\right),$$
$$B(t) = b\left(\frac{T}{2}(t+1)\right), \quad F(t) = f\left(\frac{T}{2}(t+1)\right), \tag{14}$$

then, Eq. (13) can be reduced to

$$Y(t) = \left(\frac{T}{2}\right)^\alpha \frac{A(t)}{\Gamma(\alpha)} \int_{-1}^t B(\mu)(t-\mu)^{\alpha-1} Y(\mu)\, d\mu + F(t), \quad t \in \Lambda = [-1,1]. \tag{15}$$

At last, under the following linear transformation

$$\mu = \mu(t,\theta) = \frac{t+1}{2}\theta + \frac{t-1}{2}, \quad \theta \in \Lambda, \tag{16}$$

the Eq. (15) becomes

$$Y(t) = \left(\frac{T}{2}\right)^\alpha \left(\frac{t+1}{2}\right)^\alpha \frac{A(t)}{\Gamma(\alpha)} \int_{-1}^{1} B(\mu(t,\theta)) (1-\theta)^{\alpha-1} Y(\mu(t,\theta)) d\theta + F(t). \tag{17}$$

The Legendre spectral collocation method for Eq. (17) is to seek $U \in S_N(\Lambda)$ ($N \geq 1$), such that

$$U(t) = \frac{T^\alpha}{4^\alpha \Gamma(\alpha)} p_{t,N} \left\{ (t+1)^\alpha A(t) \int_{-1}^{1} (1-\theta)^{\alpha-1} p_{\theta,N}^{\alpha-1,0} \left( B(\mu(t,\theta)) U(\mu(t,\theta)) \right) d\theta \right\}$$
$$+ p_{t,N}(F(t)). \tag{18}$$

Now, we want to describe a numerical implementation for Eq. (18). To this end, we use the following relations

$$U(t) = \sum_{i=0}^{N} u_i L_i(t),$$

$$p_{t,N} \left( p_{\theta,N}^{\alpha-1,0} \left( (t+1)^\alpha A(t) B(\mu(t,\theta)) U(\mu(t,\theta)) \right) \right) \sum_{i=0}^{N} \sum_{j=0}^{N} d_{ij} L_i(t) J_j^{\alpha-1,0}(\theta). \tag{19}$$

Then by using (19) and then Eq. (3), a direct computation leads to the following equation as

$$\frac{T^\alpha}{4^\alpha \Gamma(\alpha)} \int_{-1}^{1} (1-\theta)^{\alpha-1} p_{x,N} \left( p_{\theta,N}^{\alpha-1,0} \left( (t+1)^\alpha A(t) B(\mu(t,\theta)) U(\mu(t,\theta)) \right) \right) d\theta$$
$$= \frac{T^\alpha}{4^\alpha \Gamma(\alpha)} \sum_{i=0}^{N} \sum_{j=0}^{N} d_{ij} L_i(t) \int_{-1}^{1} (1-\theta)^{\alpha-1} J_j^{\alpha-1,0}(\theta) d\theta$$
$$= \frac{T^\alpha}{2^\alpha \Gamma(\alpha+1)} \sum_{i=0}^{N} d_{i0} L_i(t). \tag{20}$$

Applying Eqs. (3) to (6), $d_{i0}$ can be readily obtained as

$$d_{i0} = \frac{\alpha(2i+1)}{2^{\alpha+1}} \sum_{l_1}^{N} \sum_{l_2}^{N} (x_{l_1} + 1)^\alpha A(x_{l_1}) B\left(\mu(x_{l_1}, \theta_{l_2}^{\alpha-1,0})\right) U\left(\mu(x_{l_1}, \theta_{l_2}^{\alpha-1,0})\right) L_i(x_{l_1}) \omega_{l_1} \omega_{l_2}^{\alpha-1,0}.$$

Hence, by inserting Eqs. (19) and (20) in Eq. (18), we conclude

$$\sum_{i=0}^{N} u_i L_i(t) = \frac{T^\alpha}{2^\alpha \Gamma(\alpha+1)} \sum_{i=0}^{N} d_{i0} L_i(t). \tag{21}$$

Finally, by comparing the expansion coefficients of eq. (21), we can get

$$u_i = \frac{T^\alpha}{2^\alpha \, \Gamma(\alpha+1)} d_{i0}, \qquad 0 \leq i \leq N. \tag{22}$$

By replacing the above relation, we obtain the expression of $U(t)$ accordingly.

## 4. Convergence Analysis

In this section, we want to carry out the error analysis for the numerical method (18) under $L^2(\Lambda)$ and $L^\infty(\Lambda)$. Here, we present some preparations as follow [14,17]:

$L^2_{\omega^{q_1,q_2}}(\Lambda)$: the Jacobi-weighted $L^2$ Hilbert space with the scalar product

$$(u,v) = \int_{-1}^{1} u(t)\, v(t)\, \omega(t) dt, \quad \forall u, v \in L^2_{\omega^{q_1,q_2}}(\Lambda).$$

Therefore, the corresponding norm is

$$\|u\|_{L^2_{\omega^{q_1,q_2}}} = (u,u)^{\frac{1}{2}}.$$

$H^m_{\omega^{q_1,q_2}}(\Lambda) = \left\{ u \in L^2_{\omega^{q_1,q_2}}(\Lambda) : \frac{d^i u}{dx^i} \in L^2_{\omega^{q_1,q_2}}(\Lambda),\ i = 0,1,\dots m \right\}$: the Jacobi-weighted Sobolev spaces with the norm and semi norm respectively

$$\|u\|^2_{H^m_{\omega^{q_1,q_2}}} = \sum_{k=0}^{m} \|\partial^k_t u\|^2_{L^2_{\omega^{q_1,q_2}}},$$

and

$$|u|_{H^m_{\omega^{q_1,q_2}}} = \|\partial^k_t u\|_{\omega^{q_1+k,q_2+k}}$$

$H^m_{\omega^{q_1,q_2}}$ by its inner product is Hilbert space. For the purpose of convenience, we write $L^2(\Lambda) = H^0_{\omega^{0,0}}$, $\|\cdot\|_2 = \|\cdot\|_{L^2(\Lambda)}$ and $\|\cdot\|_\infty = \|\cdot\|_{L^\infty(\Lambda)}$.

**Lemma 4.1 [17]** Assume that $u \in H^m_{\omega^{q_1,q_2}}$ with $q_1, q_2 > -1$, and $m \geq 1$. Then for any integer number $k$ such that $0 \leq k \leq m \leq N+1$, exists the constant $C > 0$ that the following relation holds:

$$\left\| \partial^k_t \left( u - p^{q_1,q_2}_{t,N}(u) \right) \right\|_{\omega^{q_1+k,q_2+k}} \leq CN^{k-m}\, \|\partial^m_t u\|_{\omega^{q_1+k,q_2+k}}. \tag{23}$$

In particular, for any u ∈ $H^m(\Lambda)$ with $1 \leq m \leq N+1$, the above relation reduced to the following result

$$\|u - p_{t,N}(u)\|_{H^1(\Lambda)} \leq C N^{\frac{3}{2}-m} \|\partial_t^m u\|_2. \tag{24}$$

**Theorem 4.2** For sufficiently large $N$, the Legendre spectral-collocation approximation converges to exact solution in $L^2$-norm, i.e.

$$\|e_N\|_2 = \|Y - U\|_2 \to 0.$$

Proof. Assume that $U(x)$ is obtained by using the Legendre spectral-collocation method Eq. (1). Then, according to definition of $p_{t,N}$ in last section, we will have

$$\|e_N\|_2 \leq \|Y - p_{t,N}(Y)\|_2 + \|p_{t,N}(Y) - U\|_2. \tag{25}$$

By using Lemma 4.1, we can get for any integer $1 \leq m \leq N+1$,

$$\|Y - p_{t,N}(Y)\|_2 \leq C N^{-m} \|\partial_t^m Y\|_{\omega^{m,m}}. \tag{26}$$

The rest of our proof is to show $\|p_{t,N}(Y) - U\|_2 \to 0$. To this end, along with definition $Y$ in Eq. (17), for $N \geq 1$ we have

$$p_{t,N} Y(t) = \frac{T^\alpha}{2^\alpha \Gamma(\alpha)} p_{t,N}\left(A(t) \int_{-1}^{t} (t-\mu)^\alpha F(B(\mu), Y(\mu)) d\mu\right) + p_{t,N}(F), \tag{27}$$

where, $F(B(\mu), Y(\mu)) = B(\mu(t,\theta)) Y(\mu(t,\theta))$. Similarly, we can write

$$U(t) = \frac{T^\alpha}{2^\alpha \Gamma(\alpha)} p_{t,N}\left\{A(t) \int_{-1}^{t} (t-\mu)^{\alpha-1} p_{\mu,N}^{\alpha-1,0}\left(F(B(\mu), U(\mu))\right) d\mu\right\} + p_{t,N}(F(t)), \tag{28}$$

where

$$p_{\mu,N}^{\alpha-1,0} v(\mu) = p_{\theta,N}^{\alpha-1,0} v(\mu(t,\theta))\Big|_{\theta = \frac{2\mu}{t+1} - \frac{t-1}{t+1}}. \tag{29}$$

Before we proceed to the rest of proof, we present some relations that will be used later. By virtue of above definition and then using standard Jacobi-Gauss quadrature formula (4), we obtain

$$\int_{-1}^{t} (t-\mu)^{\alpha-1} p_{\mu,N}^{\alpha-1,0} v(\mu) d\mu = \left(\frac{1+t}{2}\right)^\alpha \int_{-1}^{1} (1-\theta)^{\alpha-1} p_{\theta,N}^{\alpha-1,0} v(\mu(t,\theta)) d\theta$$

$$= \left(\frac{1+t}{2}\right)^\alpha \sum_{j=0}^{N} v\left(\mu(t, \theta_j^{\alpha-1,0})\right) \omega_j^{\alpha-1,0}$$

$$= \left(\frac{1+t}{2}\right)^\alpha \sum_{j=0}^{N} v(\mu_j^{\alpha-1,0}) \omega_j^{\alpha-1,0}. \tag{30}$$

Again, in the same scheme we have

$$\int_{-1}^{t}(t-\mu)^{\alpha-1}\left(p_{\mu,N}^{\alpha-1,0}v(\mu)\right)^2 d\mu = \left(\frac{1+t}{2}\right)^{\alpha}\sum_{j=0}^{N}v^2(\mu_j^{\alpha-1,0})\omega_j^{\alpha-1,0}. \tag{31}$$

At last, by applying Lemma 4.1 and Eq. (29) together, for any integer $m$ ($1 \leq m \leq N+1$), we conclude that

$$\int_{-1}^{t}(t-\mu)^{\alpha-1}\left(v(\mu) - p_{\mu,N}^{\alpha-1,0}v(\mu)\right)^2 d\mu$$

$$= \left(\frac{1+t}{2}\right)^{\alpha}\int_{-1}^{1}(1-\theta)^{\alpha-1}\left(v(\mu(t,\theta)) - p_{\theta,N}^{\alpha-1,0}v(\mu(t,\theta))\right)^2 d\theta$$

$$\leq CN^{-m}\left(\frac{1+t}{2}\right)^{\alpha}\int_{-1}^{1}\left(\partial_\theta^m v(\mu(t,\theta))\right)^2 (1-\theta)^{\alpha+m-1}(1+\theta)^m d\theta$$

$$= CN^{-m}\int_{-1}^{x}\left(\partial_\mu^m v(\mu)\right)^2 (x-\mu)^{\alpha+m-1}(1+\mu)^m d\mu. \tag{32}$$

Now, let us explain the rest of proof. By subtracting Eq. (28) from (27), we derive that

$$p_{t,N}Y(t) - U(t)$$

$$= \frac{T^\alpha}{2^\alpha \Gamma(\alpha)} \times p_{t,N}\left(A(t)\int_{-1}^{t}(t-\mu)^{\alpha-1}\left(F(B(\mu),Y(\mu)) - p_{\mu,N}^{\alpha-1,0}\left(F(B(\mu),U(\mu))\right)\right) d\mu\right).$$

The above formula can be rewritten as

$$p_{t,N}Y(t) - U(t)$$

$$= \frac{T^\alpha}{2^\alpha \Gamma(\alpha)} \times p_{t,N}\left(A(t)\int_{-1}^{t}(t-\mu)^{\alpha-1}\left(F(B(\mu),Y(\mu)) - p_{\mu,N}^{\alpha-1,0}\left(F(B(\mu),Y(\mu))\right)\right) d\mu\right)$$

$$+ \frac{T^\alpha}{2^\alpha \Gamma(\alpha)}$$

$$\times p_{t,N}\left(A(t)\int_{-1}^{t}(t-\mu)^{\alpha-1}p_{\mu,N}^{\alpha-1,0}\left(F(B(\mu),Y(\mu)) - F(B(\mu),U(\mu))\right) d\mu\right). \tag{33}$$

We let

$$E_1(t) = \frac{T^\alpha}{2^\alpha \Gamma(\alpha)} \times p_{t,N}\left(A(t)\int_{-1}^{t}(t-\mu)^{\alpha-1}\left(F(B(\mu),Y(\mu)) - p_{\mu,N}^{\alpha-1,0}\left(F(B(\mu),Y(\mu))\right)\right) d\mu\right),$$

$$E_2(t) = \frac{T^\alpha}{2^\alpha \Gamma(\alpha)} \times p_{t,N}\left(A(t)\int_{-1}^{t}(t-\mu)^{\alpha-1}p_{\mu,N}^{\alpha-1,0}\left(F(B(\mu),Y(\mu)) - F(B(\mu),U(\mu))\right) d\mu\right).$$

Next, by using standard Legendre-Gauss quadrature formula we obtain

$\|E_1\|_2$

$$= \frac{T^\alpha}{2^\alpha \Gamma(\alpha)} \left\{ \sum_{j=0}^{N} \omega_j\, A(x_j) \left( \int_{-1}^{x_j} (x_j - \mu)^{\alpha-1} \left( F(B(\mu), Y(\mu)) - p_{\mu,N}^{\alpha-1,0}\big|_{x_j} \left( F(B(\mu), Y(\mu)) \right) \right) d\mu \right)^2 \right\}^{\frac{1}{2}}.$$

By using the Cauchy-Schwartz inequality and Eq. (32), we get

$\|E_1\|_2$

$$\leq \max_{t \in \Lambda} A(t)\, \frac{T^\alpha}{2^\alpha \Gamma(\alpha)} \left\{ \sum_{j=0}^{N} \omega_j \int_{-1}^{x_j} (x_j - \mu)^{\alpha-1} d\mu \int_{-1}^{x_j} (x_j - \mu)^{\alpha-1} \left( F(B(\mu), Y(\mu)) \right.\right.$$

$$\left.\left. - p_{\mu,N}^{\alpha-1,0}\big|_{x_j} \left( F(B(\mu), Y(\mu)) \right) \right)^2 d\mu \right\}^{\frac{1}{2}}$$

$$= \max_{t \in \Lambda} A(t)\, \frac{T^\alpha}{2^\alpha \Gamma(\alpha+1)}$$

$$\times \left\{ \sum_{j=0}^{N} \omega_j (x_j + 1)^\alpha \int_{-1}^{x_j} (x_j - \mu)^{\alpha-1} \left( F(B(\mu), Y(\mu)) - p_{\mu,N}^{\alpha-1,0}\big|_{x_j} \left( F(B(\mu), Y(\mu)) \right) \right)^2 d\mu \right\}^{\frac{1}{2}}$$

$$\leq C\, \max_{t \in \Lambda} A(t)\, \frac{T^\alpha}{2^\alpha \Gamma(\alpha+1)}\, N^{-m}$$

$$\times \left\{ \sum_{j=0}^{N} \omega_j (x_j + 1)^\alpha \int_{-1}^{x_j} \left( \partial_\mu^m F(B(\mu), Y(\mu)) \right)^2 (x_j - \mu)^{\alpha+m-1} (1+\mu)^m d\mu \right\}^{\frac{1}{2}}$$

$$= C_1\, N^{-m}\, \left\| \partial_\mu^m F(B(.), Y(.)) \right\|_{\omega^{\alpha+m-1,m}}. \tag{34}$$

Similarly, by applying Eqs. (4) and (31) and the Cauchy-Schwartz inequality, we can get

$\|E_2\|_2$

$$= \frac{T^\alpha}{2^\alpha \Gamma(\alpha)} \left\{ \sum_{j=0}^{N} \omega_j \left[ A(x_j) \int_{-1}^{x_j} (x_j - \mu)^{\alpha-1} p_{\mu,N}^{\alpha-1,0}\big|_{x_j} \left( F(B(\mu), Y(\mu)) - F(B(\mu), U(\mu)) \right) d\mu \right]^2 \right\}^{\frac{1}{2}}$$

$$\leq \max_{t \in \Lambda} A(t)\, \frac{T^\alpha}{2^\alpha \Gamma(\alpha)} \left\{ \sum_{j=0}^{N} \omega_j \int_{-1}^{x_j} (x_j - \mu)^{\alpha-1} d\mu \right.$$

$$\left. \times \int_{-1}^{x_j} (x_j - \mu)^{\alpha-1} \left( p_{\mu,N}^{\alpha-1,0}\big|_{x_j} \left( F(B(\mu), Y(\mu)) - F(B(\mu), U(\mu)) \right) \right)^2 d\mu \right\}^{\frac{1}{2}}$$

$$= \max_{t\in\Lambda} A(t) \frac{T^\alpha}{2^\alpha \Gamma(\alpha)} \left\{ \sum_{j=0}^{N} \omega_j \frac{(x_j+1)^\alpha}{\alpha} \sum_{i=0}^{N} \omega_i^{\alpha-1,0} \left( F\left(B(\mu_i^{\alpha-1,0}), Y(\mu_i^{\alpha-1,0})\right) \right.\right.$$

$$\left.\left. - F\left(B(\mu_i^{\alpha-1,0}), U(\mu_i^{\alpha-1,0})\right) \right)^2 \right\}^{\frac{1}{2}}. \tag{35}$$

On the other hand, Jenson's inequality for any convex function $g$ on interval $(a,b)$ is

$$g(\lambda t_1 + (1-\lambda)t_2) \leq \lambda g(t_1) + (1-\lambda)g(t_2), \quad \forall t_1, t_2 \in (a,b), \quad \forall \lambda \in [0,1].$$

Thus, since

$$\frac{d^2}{dt^2}(x_j+1)^t = (x_j+1)^t \ln^2(x_j+1) > 0, \quad x_j \text{ is constant,}$$

then, for any given $x_j \in (-1,1)$, $g(t) = (x_j+1)^t$ is convex function and according to Jenson's inequality on interval $[0,1]$, we will have

$$g(t) = g((1-t)\times 0 + t\times 1) \leq (1-t)g(0) + tg(1). \tag{36}$$

Therefore, by applying the previous inequality, we yield the following inequality

$$\sum_{j=0}^{N} \omega_j (x_j+1)^\alpha \leq \sum_{j=0}^{N} \omega_j \{(1-\alpha) + \alpha(x_j+1)\}$$

$$= (1-\alpha) \sum_{j=0}^{N} \omega_j + \alpha \sum_{j=0}^{N} \omega_j(x_j+1)$$

$$= (1-\alpha) \int_\Lambda dt + \alpha \int_\Lambda (x+1) dt$$

$$= 2, \quad \forall \alpha \in [0,1]. \tag{37}$$

Next, by using relation (37), Eq. (35) will be reduced to

$$\|E_2\|_2$$

$$\leq \max_{t\in\Lambda} A(t) \frac{T^\alpha}{2^{\alpha-1} \Gamma(\alpha+1)}$$

$$\times \max_{0\leq j\leq N} \left\{ \int_{-1}^{x_j} (x_j-\mu)^{\alpha-1} \left( p_{\mu,N}^{\alpha-1,0}\big|_{x_j} \left(F(B(\mu),Y(\mu)) - F(B(\mu),U(\mu))\right) \right)^2 d\mu \right\}^{\frac{1}{2}}$$

$$\leq \max_{t\in\Lambda} A(t) \max_{t\in\Lambda} B(t) \frac{T^\alpha}{2^{\alpha-1} \Gamma(\alpha+1)}$$

$$\times \max_{0 \le j \le N} \left\{ \int_{-1}^{x_j} (x_j - \mu)^{\alpha-1} \left( p_{\mu,N}^{\alpha-1,0} \big|_{x_j} (Y(\mu) - U(\mu)) \right)^2 d\mu \right\}^{\frac{1}{2}}$$

$$\le \max_{t \in \Lambda} A(t) \max_{t \in \Lambda} B(t) \frac{T^\alpha}{2^{\alpha-1} \Gamma(\alpha+1)}$$

$$\times \max_{0 \le j \le N} \left\{ \left( \int_{-1}^{x_j} (x_j - \mu)^{\alpha-1} \left( p_{\mu,N}^{\alpha-1,0} \big|_{x_j} (Y(\mu)) - Y(\mu) \right)^2 d\mu \right)^{\frac{1}{2}} \right.$$

$$\left. + \left( \int_{-1}^{x_j} (x_j - \mu)^{\alpha-1} (Y(\mu) - U(\mu))^2 d\mu \right)^{\frac{1}{2}} \right\}. \tag{38}$$

At last, along with Lemma 4.1 and relations (32), (37), the previous result yields

$$\|E_2\|_2$$

$$\le C_2 \left\{ CN^{-m} \max_{0 \le j \le N} \left( \int_{-1}^{x_j} \left( \partial_\mu^m Y(\mu) \right)^2 (x_j - \mu)^{\alpha+m-1} (1+\mu)^m d\mu \right)^{\frac{1}{2}} \right.$$

$$\left. + \frac{2}{\alpha} \left( \int_{-1}^{1} (Y(\mu) - U(\mu))^2 d\mu \right)^{\frac{1}{2}} \right\}$$

$$\le C_2 C N^{-m} \|\partial_\mu^m Y\|_{\omega^{\alpha+m-1,m}} + \frac{2C_2}{\alpha} \|e_N\|. \tag{39}$$

Finally, relations $(33) - (34), (39)$ prove that the approximation is convergent in $L^2$ −norm. Hence the theorem is proved.

## 5. Numerical example

To show efficiency of our numerical method, the following examples which have exact solution are considered.

**Example 5.1.** Consider the following fractional integral equation [20]

$$y(t) = \frac{0.01}{\Gamma(0.5)} t^{\frac{5}{2}} \int_0^t (t-s)^{-\frac{1}{2}} y(s) \, ds + \sqrt{\pi}(1+t)^{-\frac{3}{2}} - 0.02 \frac{t^3}{1+t}, \quad t \in [0,1].$$

The corresponding exact solution is given by $y(t) = \sqrt{\pi} \, (1+t)^{-\frac{3}{2}}$. Figure 1 present the approximate and exact solutions which are found in very good agreement. In Figure 2, the numerical errors are plotted in $L^2$ −norm for $2 \le N \le 24$. As expected, an exponential rate of convergence is observed which confirmed our theoretical predictions in comparison of results in [20].

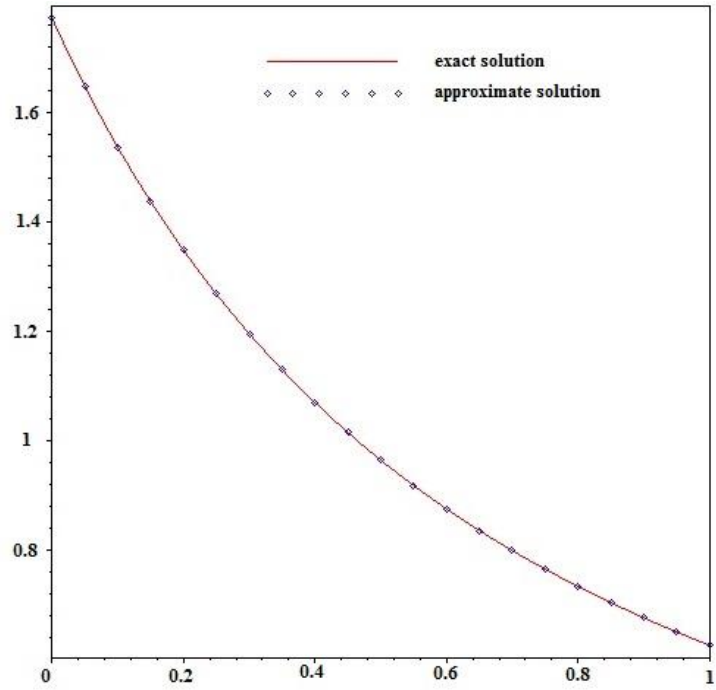

**Figure 1.** Comparison between exact solution and approximate solution of Example 5.1 for N = 10

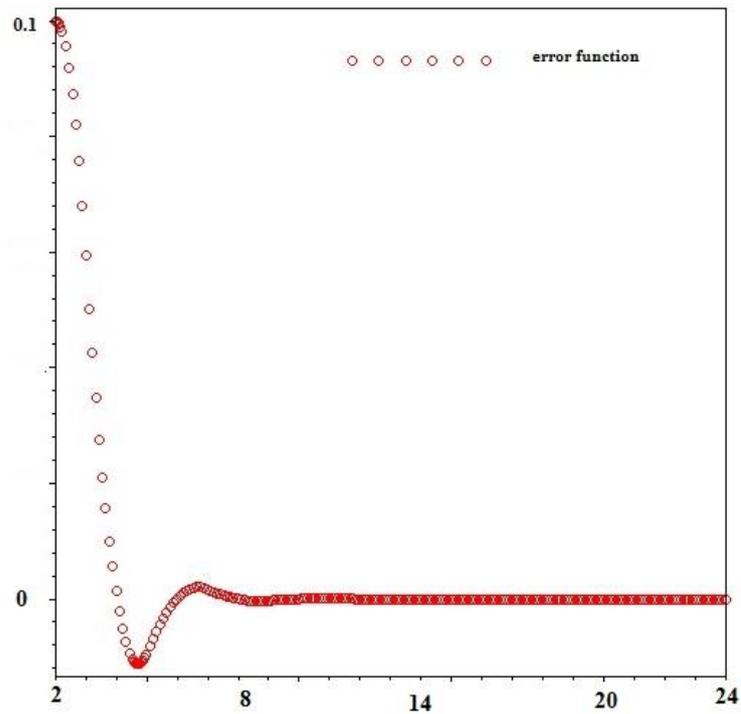

**Figure 2.** The error function of Example 5.1 versus the number of interpolation operator

**Example 5.2.** Our second example is the following fractional integral equation [20]

$$y(t) = \frac{1}{27\,\Gamma\left(\frac{2}{3}\right)} \int_0^t s(t-s)^{-\frac{1}{3}} y(s)\, ds + \Gamma\left(\frac{2}{3}\right) t - \frac{1}{40} t^{\frac{8}{3}}, \quad t \in [0,1].$$

The exact solution is $y(x) = \Gamma\left(\frac{2}{3}\right)t$. We get the numerical solution of the above fractional integral equation as

$$y_1(t) = \Gamma\left(\frac{2}{3}\right)t,$$

$$y_2(t) = 1.21 \times 10^{-2} + \Gamma\left(\frac{2}{3}\right)t,$$

$$y_4(t) = 2.15 \times 10^{-9} + \Gamma\left(\frac{2}{3}\right)t + 3.12 \times 10^{-10}\,t^2 + 2.39 \times 10^{-7}\,t^3 + 4.81 \times 10^{-9}.$$

The numerical results of proposed method can be seen from **Figure 3** and **Figure 4**. These results indicate that the spectral accuracy is obtained for this problem, although the given function $f(t)$ is not very smooth.

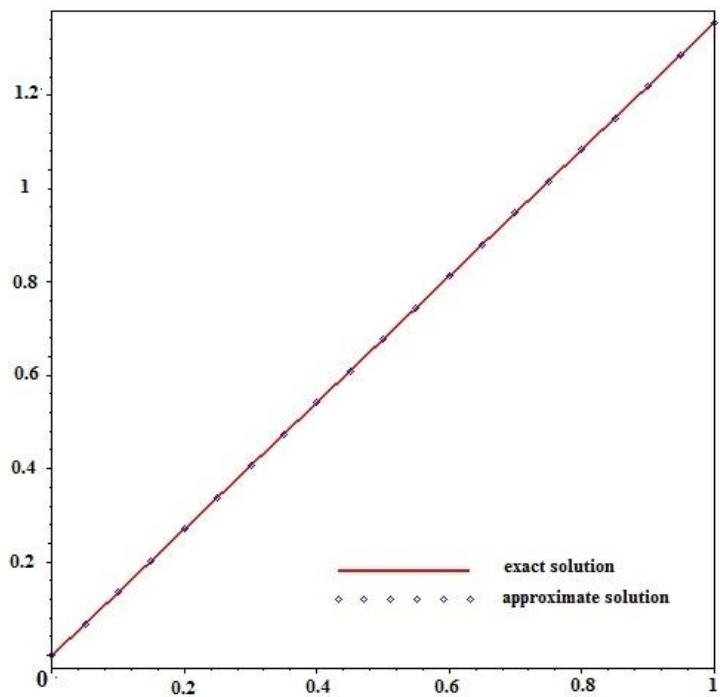

**Figure 3.** Comparison between exact solution and approximate solution of Example 5.2 for $N = 6$

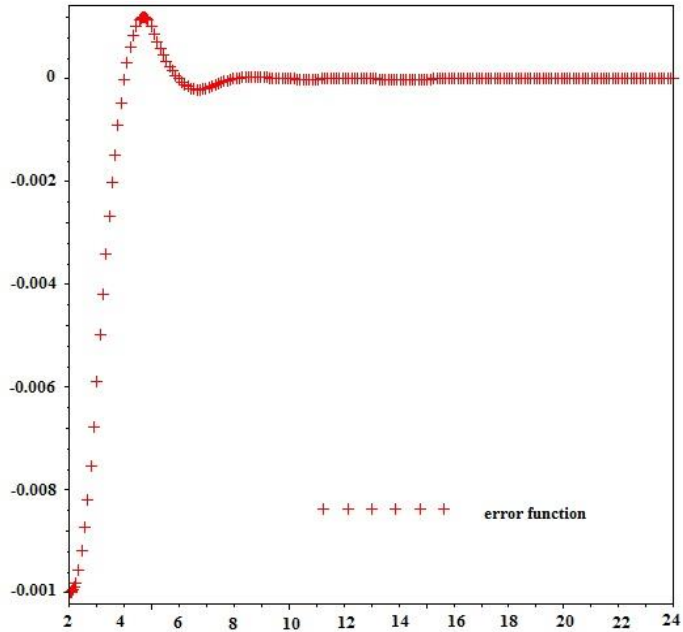

**Figure 4.** The error function of Example 5.1 versus the number of interpolation operator

These results indicate that the spectral accuracy is obtained by using Legendre-Collocation method for this problem.

**Example 5.3.** Consider the following fractional integral equation [20]

$$y(t) = -\frac{1}{\Gamma(0.5)} \int_0^t (t-s)^{-\frac{1}{2}} y(s)\, ds + 2\sqrt{\frac{t}{\pi}}, \quad t \in [0,1].$$

The exact solution is $y(t) = 1 - e^t\, erfc(\sqrt{t})$, such that $erfc$ is the complementary error function defined by

$$erfc(t) = 1 - erf(t) = \frac{2}{\sqrt{\pi}} \int_t^\infty e^{-s^2}\, ds.$$

We list the numerical results in Table 1 to be able to test the convergence rate of the suggested. In Figure 5, we can observe that our numerical solutions coincide closely with the exact ones. At a glance, we can find out the results of Legendre spectral-collocation method are satisfactory in the few steps.

| t | Proposed method at $N=8$ | Proposed method at $N=10$ | Exact solution |
|---|---|---|---|
| 0 | .0000000000 | .0000000000 | .0000000000 |
| 0.1 | .2764215215 | .2764215613 | .2764215616 |
| 0.2 | .3562117262 | .3562117272 | .3562117278 |
| 0.3 | .4079815809 | .4079815886 | .4079815886 |
| 0.4 | .4463937357 | .4463937452 | .4463937458 |
| 0.5 | .4768434165 | .4768434157 | .4768434159 |
| 0.6 | .5019754248 | .5019754315 | .5019754315 |
| 0.7 | .5232972359 | .5232972687 | .5232972687 |
| 0.8 | .5417539715 | .5417539771 | .5417539772 |
| 0.9 | .5579785802 | .5579785884 | .5579785885 |
| 1 | .5724164231 | .5724164243 | .5724164241 |

**Table 1.** Comparison of approximate solutions and exact solutions for various numbers N

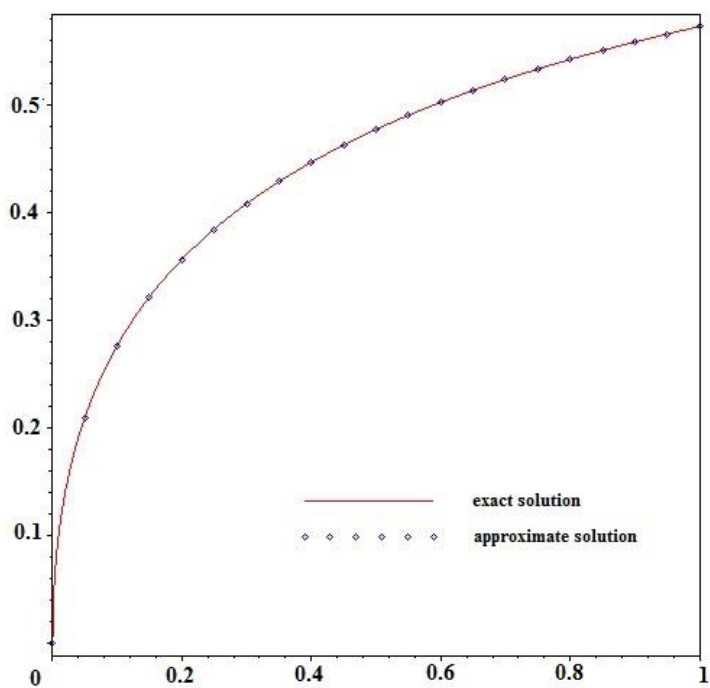

**Figure 2.** The error function of Example 5.1 versus the number of interpolation operator for $N=16$

## 6. Conclusion

In this paper we present a spectral Legendre-collocation approximation of a class of fractional integral equations of the second kind. The most important contribution of this paper is that the errors of approximations decay exponentially in $L^2 - norm$. We prove that our proposed method is effective and has high convergence rate and better than some new methods which be define in this paper. During three numerical applications, we ensure that the proposed method is efficient and powerful numerical scheme for solving many fractional problems.

**Reference**

[1] A.A. Kilbas, H.M. Srivastava and J.J. Trujillo, Theory and Applications of Fractional Differential Equations, North-Holland Math. Stud., vol. 204, Elsevier Science B.V., Amsterdam, (2006).
[2] I. Podlubny, Fractional Differential Equations, Academic Press, San Diego, (1999).
[3] A. H. Bhrawy and M. A. Alghamdi, A shifted Jacobi-Gauss-Lobatto collocation method for solving nonlinear fractional Langevin equation involving two fractional orders in different intervals, Boundary Value Problems, vol. 2012, article 62, 13 pages, 2012.
[4] E. Rawashdeh, Numerical solution of fractional integro-differential equations by collocation method, Appl. Math. Comput. 176, 1-6 (2006).
[5] Y. Nawaz, Variational iteration method and Homotopy perturbation method for fourth-order fractional integro-differential equations, Comput. Math. Appl. 61, 2330-2341 (2011).
[6] S. Momani and M. A. Noor, Numerical method for fourth-order fractional integro-differential equations, Appl. Math. Comput. 182, 754-760 (2006).
[7] H. Saeedi and F. Samimi, He's homotopy perturbation method for nonlinear Fredholm integro-differential equation of fractional order, International J. of Eng. Res. And App., Vol. 2, no. 5, pp. 52-56, 2012.
[8] A. Ahmed and S. A. H. Salh, Generalized Taylor matrix method for solving linear integro-fractional differential equations of Volterra type, App. Math. Sci., Vol. 5, no. 33-36, pp. 1765-1780, 2011.
[9] A. Arikoglo and L. Ozkol, Solution of fractional integro-differential equation by using fractional differential transform method, Chaos Solitons Fractals 34, 1473-1481 (2007).
[10] J. Wei and T. Tian, Numerical solution of nonlinear Volterra integro-differential equations of fractional order by the reproducing kernel method, Appl. Math. Model. 39, 4871-4876 (2015).
[11] Li Zhu and Qibin Fan, Solving fractional nonlinear Fredholm integro-differential equations by the second kind Chebyshev wavelet, Common. Nonlinear Sci. Numer. Simul. 17, 2333–2341 (2012).
[12] H. Saeedi, M. Mohseni Moghadam, N. Mollahasani and GN. Chuev, A CAS wavelet method for solving nonlinear Fredholm integro-differential equation of fractional order, Common. Nonlinear Sci. Numer. Simul. 16, 1154-1163 (2011).
[13] Z. Meng, L. Wang, H. Li and W. Zhang, Legendre wavelets method for solving fractional integro-differential equations, Int. J. Comput. Math. 92, 1275-1291 (2015).
[14] C. Canuto, M. Y. Hussaini, A. Quarteroni and T. A. Zang, Spectral methods: Fundamentals in single domain, Springer, Berlin (2006).
[15] P. E. Bjorstad and B. P. Tjostheim, Efficient algorithm for solving a fourth-order equation with the spectral Galerkin method, SIAM J. Sci. Comput., 18, pp. 621-632 (1997).
[16] J. Shen, Efficient spectral-Galerkin method I. Direct solvers for second- and fourth-order equations by using Legendre polynomials, SIAM J. Sci. Comput., 15, pp. 1489-1505 (1994).
[17] J. Shen, T. Tang and Li-Lian Wang, Spectral Methods: Algorithms, Analysis and Applications, Springer-Verlag, Berlin Heidelberg (2011).
[18] A. Yousefi, T. Mahdavi-Rad and S. G. Shafiei, A quadrature Tau method for solving Fractional integro-differential equations in the Caputo sense, J. of Math. Computer Sci. 16, 97-107 (2015).

[19] W. G. El-Sayed and A. M.A. El-Sayed, On the fractional integral equations of mixed type integro-differential equations of fractional orders, App. Math. Comput. 154, 461-467 (2004).

[20] S. Micula, An iterative numerical method for fractional integral equations of the second kind, J. of. Comput. and Appl. Math. (2017) DOI: http://doi.org/10.1016/j.cam.2017.12.006.

[21] Y. Yang, Y. Chen and Y. Huang, Convergence analysis of the Jacobi spectral-collocation method for fractional integro-differential equations, J. of Acta Math. Sci. 34B (3), 673-690 (2014).

[22] G. Mastroianni and D. Occorsto, Optimal systems of nodes for Lagrange interpolation on bounded intervals: a survey, J. Comput. Appli. Math. 134, 325-341 (2010).